\documentclass[12pt]{amsart}

\usepackage{amsmath,amsthm,amsfonts,amssymb,eucal}

\renewcommand {\a}{ \alpha }

\newcommand{\y}{\eta}
\newcommand{\e}{\epsilon}
\newcommand{\vare}{\varepsilon}
\newcommand{\g}{\gamma}
\newcommand{\G}{\Gamma}

\newcommand{\varf}{\varphi}
\renewcommand{\d}{\delta}
\newcommand{\D}{\Delta}
\newcommand{\s}{\sigma}
\renewcommand{\l}{\lambda}
\renewcommand{\L}{\Lambda}
\newcommand{\z}{\zeta}

\newcommand{\Om}{\Omega}

\newcommand{\R}{ \mathbb R}
\newcommand{\C}{ \mathbb C}
\newcommand{\N}{ \mathbb N}
\newcommand{\Sq}{ \mathbb S}
\newcommand{\Z}{ \mathbb Z}

\newcommand{\CA}{\mathcal A}

\newcommand{\CB}{\mathcal B}

\newcommand{\CD}{\mathcal D}

\newcommand{\CM}{\mathcal M}

\newcommand {\gm}{\mathfrak m}

\newcommand {\gd}{\mathfrak d}

\newcommand {\ba}{\mathbf a}

\newcommand {\bl}{\mathbf l}

\newcommand {\bb}{\mathbf b}
\newcommand {\bm}{\mathbf m}
\newcommand {\br}{\mathbf r}

\newcommand {\BA}{\mathbf A}

\newcommand {\BJ}{\mathbf J}

\newcommand {\BL}{\mathbf L}
\newcommand {\BM}{\mathbf M}

\newcommand {\BR}{\mathbf R}

 \DeclareMathOperator{\im}{Im}
\DeclareMathOperator{\re}{Re} 
\DeclareMathOperator{\res}{\restriction}

\newtheorem{thm}{Theorem}[section]

\newtheorem{lem}[thm]{Lemma}
\newtheorem{prop}[thm]{Proposition}

\theoremstyle{definition}
\newtheorem{defn}[thm]{Definition}%[section]

\theoremstyle{remark}
\newtheorem{rem}[thm]{Remark}

\numberwithin{equation}{section}

\newcommand{\thmref}[1]{Theorem~\ref{#1}}

\begin{document}

\title[A family of differential operators]{On a family of differential
operators with the coupling parameter in the boundary
condition}
 \dedicatory{Dedicated to Professor Des Evans on the occasion of
his 65-th
 birthday}
\author[Rozenblum]{G. Rozenblum}
\address[G. Rozenblum]{Department of Mathematics \\
                        Chalmers University of Technology
                        and University of Gothenburg \\
                         S-412 96, Gothenburg,
                        Sweden}
\email{grigori@math.chalmers.se}
\author[Solomyak]{M. Solomyak}
\address[M. Solomyak]{Department of Theoretical Mathematics\\ The Weizmann Institute of Science\\
Rehovot, 76100, Israel} \email{michail.solomyak@weizmann.ac.il}
\begin{abstract}

We study a family of differential operators $\BL_\a$
in two variables, depending on the coupling parameter
$\a\ge0$ that appears only in the boundary conditions.
Our main concern is the spectral properties of
$\BL_\a$, which turn out to be  quite different for
$\a<1$ and for $\a>1$. In particular,
 $\BL_\a$ has a unique self-adjoint realization for $\a<1$ and
many such realizations for $\a>1$. In the more
difficult case $\a>1$ an analysis of non-elliptic
pseudodifferential operators in dimension one is
involved.
\end{abstract}

\maketitle

\section{Introduction}\label{int}
In the paper \cite{SM} Smilansky suggested  a mathematical model
which he called "The irreversible quantum graph".  In this model a
one-dimensional quantum graph interacts with a finite system of
harmonic oscillators attached at different points of the graph.
Regardless of the physical meaning of this model, it is quite
interesting from the mathematical point of view, since, being a
singular perturbation problem, it exhibits many unusual effects.
These effects appear already in the one-oscillator case. They were
discussed in the survey paper \cite{S3}, see also references
therein.

In the simplest case (the graph is a real line, with
only one oscillator attached) the problem consists in
the study of a family of differential operators
$\BA_\a$ on $\R^2$, depending on the coupling
parameter $\a\ge0$. The differential expression which
defines the action of $\BA_\a$ does not involve $\a$,
this parameter appears only  in the transmission
condition across a straight line in the plane. The
operator $\BA_0$ admits an exhaustive description via
the separation of variables, and the passage to
$\BA_\a$ with $\a \ne0$ can be expressed, at least
formally, in the terms of perturbations of quadratic
forms. The main peculiarity of the problem stems from
the fact that the perturbation is too strong: it is
only relatively bounded but not relatively compact
with respect to the operator $\BA_0$ (in the sense of
quadratic forms). For this reason, the standard
machinery of the perturbation theory does not work.
Still, it turned out to be possible to give a detailed
description of the spectrum $\s(\BA_\a)$ for all
$\a>0$. A borderline value $\a^*$ of the parameter
$\a$ exists, such that the properties of $\s(\BA_\a)$
are quite different for $\a<\a^*$ and for $\a\ge\a^*$.
For $\a<\a^*$ the absolutely continuous (a.c.)
spectrum of $\BA_\a$ is the same as for $\BA_0$,
including the multiplicity. Eigenvalues appear below
the bottom of $\s(\BA_0)$, their number grows
indefinitely as $\a\nearrow\a^*$ and satisfies an
asymptotic relation of a non-standard type. For
$\a=\a^*$ these eigenvalues disappear and a new branch
of the a.c. spectrum appears instead, filling
$[0,\infty)$. For $\a$ above the threshold $\a^*$, the
operator $\BA_\a$ is not semi-bounded any more and its
a.c. spectrum fills the whole real line. Thus, the
system exhibits a sort of phase transition as the
parameter $\a$ crosses the threshold $\a^*$.

The mathematical mechanism behind such a behaviour of the spectrum
lies in a very special form of the transmission condition for the
operator $\BA_\a$. This condition generates in a natural way an
infinite Jacobi matrix which depends on the parameter $\a$ and
whose spectral properties for $\a<\a^*$ and for $\a\ge\a^*$ are
quite different.

The papers \cite{ESI} and \cite{ESII} are devoted to the case of
two oscillators, but actually their results show what happens in
the general case of an arbitrary number of oscillators. It was an initiative
of Des Evans, to start the work on these papers, and we take pleasure
in emphasizing his role in the study of this class of problems.

In the present paper we investigate another family of differential
operators, say $\BL_\a$, of a similar nature. It was also proposed
by Smilansky (private communication). Again, all operators in the
family are determined by a differential expression not depending
on the parameter, and they differ by the transmission condition.
Like in the case of the family $\BA_\a$, a certain family of
Jacobi matrices is closely related to the operator. However, the
properties of the two families are rather different and another
type of phase transition occurs. Namely, for large values of $\a$
the operator $\BL_\a$ has many self-adjoint realizations, and the
negative spectrum of each realization is discrete and unbounded
from below. The mechanism of this transition lies in an unusual
breaking of the Shapiro -- Lopatinsky ellipticity condition in
several points on the interface line, and the analysis of this
situation involves a study of {\it a priori} estimates for some
non-elliptic pseudodifferential operators.

In the last section of the paper we briefly consider
yet another family $\BM_\a$ of differential operators.
It looks rather similar to the family $\BL_\a$, but
some important details in the behaviour of the
spectrum are quite different.

Taken together, the families $\BA_\a$, $\BL_\a$, and
$\BM_\a$ show that presence of the coupling parameter
in the boundary condition may cause quite different
types of the phase transition. It is tempting to
develop a general scheme which would include all these
examples as special cases.

\section{Stating the problem. Preliminaries}\label{red}

 We study a
 family $\BL_\a$ of differential operators on the
cylinder $\Om=\R\times\Sq^1$ identified with the strip
$\R\times(0, 2\pi)$ with periodic boundary conditions
for all functions involved. Further on, $x$  stands
for the co-ordinate on $\R$ and $y$ for the
co-ordinate on $\Sq^1$.  The operator $\BL_\a$ is
generated by the Laplacian $-\D
U=-U''_{x^2}-U''_{y^2}$ and two conditions at $x=0$.
The first condition is the continuity
\begin{equation}\label{1.-1}
U(0+,y)=U(0-,y)\ \left(=U(0,y)\right)
\end{equation}
and the second one is a `transmission condition' at $x=0$:
\begin{equation}\label{1.0}
U'_x(0+,y)-U'_x(0-,y)=i\a\left(U'_y(0,y)\cos y +(U(0,y)\cos
y)'_y\right).
\end{equation}
In \eqref{1.0} $\a$ is a real parameter. The passage
$\a\mapsto -\a$ corresponds to the change of variables
$y\mapsto y+\pi$, which does not affect the spectrum.
For this reason it is enough to consider $\a\ge0$.

By using the Fourier expansion
\begin{equation}\label{1.u}
U=(2\pi)^{-1/2}\sum_{n\in\Z} u_n(x)e^{in y}
\end{equation}
(in short, $U\sim\{u_n\}$), we reduce  the problem
formally to an infinite system of ordinary
differential operators on the real axis,
\begin{equation}\label{1.1}
-\D U\sim\{-u_n''+n^2u_n\},\;x\ne 0,\qquad n\in\Z,
\end{equation}
coupled by the conditions
\begin{gather}
u_n(0+)=u_n(0-)\ \left(=u_n(0)\right),\label{1.1x}\\
u_n'(0+)-u_n'(0-)=-\a\bigl((n+1/2)u_{n+1}(0)+(n-1/2)u_{n-1}(0)\bigr).
\label{1.2}
\end{gather}

The operator $\BL_0$ is just the standard Laplacian on the
cylinder $\Om$, with the domain   $H^2(\Omega)$. Thus, for $\a=0$
the above formal reduction of the partial differential operator is
legal, the system decouples, and we get

\begin{equation}\label{1.l0}
    \BL_0={\sum_{n\in\Z}}^\oplus\bigl(-\frac{d^2}{dx^2}+n^2\bigr).
\end{equation}
Here $-\frac{d^2}{dx^2}$ stands for the self-adjoint
operator in $L^2(\R)$ with the domain $H^2(\R)$ and
the symbol $\sum^\oplus$ denotes the orthogonal sum of
operators. The expansion \eqref{1.l0} leads to the
complete description of the spectrum $\s(\BL_0)$: it
is absolutely continuous, fills the half-line
$[0,\infty)$, and its multiplicity function is given
by
\begin{equation}\label{1.ac}
 \gm_{a.c.}(\l;\BL_0)=2+4[\l], \qquad\forall \l\ge 0,
\end{equation}
where, as usual, $[\l]$ denotes the integer part of a real number
$\l$.

For $\a\ne0$ we must first  specify in what sense the
conditions \eqref{1.-1} and \eqref{1.0} are understood. Suppose
that $U\in L^2(\Omega)$ is a weak solution of the equation $-\D
U=F\in L^2$ in each semi-cylinder \[\Omega_\pm=\{(x,y)\in\Om:\pm
x>0\}.\]
 Take any $\L\in\C\setminus[0,\infty)$ and consider the function
 \begin{equation}\label{1.u0}
    U_0=(\BL_0-\L)^{-1}(F-\L U)\in H^2(\Omega).
\end{equation} The function
$W=U-U_0$ belongs to $L^2(\Om)$ and satisfies the equation $\D
W+\L W=0$ in each semi-cylinder $\Omega_\pm$. Let $W_\pm$ stand
for the restriction of $W$ to  $\Omega_\pm$. The functions $W_\pm$
can be expanded in the Fourier series
\begin{equation}\label{1A:Fourier}
    W_\pm(x,y)=\sum_{n} w_n^\pm e^{iny}e^{-|x|\sqrt{n^2-\L}},
    \qquad \re\sqrt{n^2-\L}\ge0.
\end{equation}
Both series series converge in $L^2(\Om_\pm)$, and
\[
\int_{\Omega_\pm}|W_\pm(x,y)|^2dxdy=
\sum_n\frac{|w_n^\pm|^2}{2\re{\sqrt{n^2-\L}}}.\] Hence, $W_\pm\in
L^2(\Om_\pm)$ is equivalent to $\sum_{n}
|w_n^\pm|^2(n^2+1)^{-\frac12}<\infty$. It follows that for each
$x\in\R$ the series in \eqref{1A:Fourier} converge in
$H^{-\frac12}(\Sq^1)$ and moreover, $W_\pm(x,\cdot)$ are
continuous as functions of $x$ with values in
$H^{-\frac12}(\Sq^1)$. The same is true for the function $U$, and
this explains the meaning of the condition \eqref{1.-1}: namely,
\begin{equation}\label{1.cont}
    U(0+,y)=U(0-,y) \ {\text{as distributions in }}\
    H^{-\frac12}(\Sq^1).
\end{equation}
 \vskip0.2cm

Denote by $\CM(\Om)$ the class of all functions $U\in L^2(\Om)$
which meet the following conditions.

\noindent 1. The distributions $\D (U\res\Om_\pm)$ are
functions in $L^2(\Om_\pm)$.

\noindent 2. The condition \eqref{1.cont} is satisfied.

For any $\L\notin[0,\infty)$ we also set
\begin{equation}
\CM_\L(\Om)=\bigl\{W\in\CM(\Om):\D W+\L W=0\ {\text{in }} \Om_\pm\bigr\}.\nonumber
\end{equation}
The Fourier expansion of any function
$W\in\CM_\L(\Om)$ has the form
\begin{equation}\label{1B:Fourier}
    W(x,y)=\sum_{n} w_n e^{iny}e^{-|x|\sqrt{n^2-\L}},
\end{equation}
that is, for the coefficients in \eqref{1A:Fourier} we
have $w_n^+=w_n^-\; (=w_n)$ and thus the function
$W(x,\cdot)$ is even in $x$. We also conclude that
\begin{equation}
W\in\CM_\L(\Om)\Longleftrightarrow\ W(0,\cdot)\in H^{-\frac12}(\Sq^1).\nonumber
\end{equation}
Let us recall that in the  terms of the Fourier
coefficients $w_n$ the latter inclusion is equivalent
to
\begin{equation}\label{1.sob}
 \sum_n|w_n|^2(n^2+1)^{-\frac12}<\infty.
\end{equation}
Differentiation in \eqref{1B:Fourier}  shows that for
any $W\in\CM_\L(\Om)$ the derivatives
$W'_y(x,\cdot),\;W'_x(x,\cdot)$ take values in the
space $H^{-\frac32}(\Sq^1)$. The first of them, being
an even function, is continuous in the topology of
this space for all $x\in\R^1$. The second one is
continuous in the topology of $H^{-\frac32}(\Sq^1)$
for $x\ge0$ and for $x\le0$ separately, and its jump
across the circle  $\{x=0\}$ is well defined as an
element in $H^{-\frac32}(\Sq^1)$. The decomposition
$U=U_0+W$, where $U_0$ is defined by \eqref{1.u0},
shows that the same is true for any $U\in\CM(\Om)$. In
particular, this gives the precise meaning to both
sides in \eqref{1.0} as distributions in
$H^{-\frac32}(\Sq^1)$.

Substituting the Fourier expansion \eqref{1B:Fourier} and its
differentiated forms into \eqref{1.0}, we arrive at the   system
\eqref{1.1}, \eqref{1.1x}, \eqref{1.2} which is equivalent to the
initial problem. \vskip0.2cm

The following version of the Green formula is implied by the above
argument.
\begin{lem}\label{green} For any $U\in \CM(\Om)$ and $V\in H^2(\Om)$
(so that $V(0,\cdot)\in H^{\frac32}(\Sq^1)$) we have
\begin{equation}\label{1.gr}\begin{split}
\left(\int_{\Om_+}+\int_{\Om_-}\right)(\D U \overline
V-U\overline{\D V})&dxdy\\=
-\int_{\Sq^1}\left(U'_x(0+,y)-U'_x(0-,y)\right)&\overline{V(0,y)}dy,
\end{split}\end{equation}
where the integrals on the left-hand side are
understood in the sense of distributions on $\Om_\pm$
and the integral on the right-hand side is understood
in the sense of distributions on $\Sq^1$.
\end{lem}

Denote by $\CB$ the differential operator appearing in
the condition \eqref{1.0}:
\begin{equation}\label{1.b}
    \CB u=i(u'_y \cos y+(u\cos y)'_y).
\end{equation}
The operator $\CB$ is symmetric as acting in the space
$L^2(\Sq^1)$. The following useful equality, which is
valid for $U\in\CM(\Om)$ satisfying \eqref{1.0} with
an arbitrary $\a\ge0$ and any $ V\in H^2(\Om)$, is a
direct consequence of Lemma~\ref{green}:
\begin{equation}\label{1.gr1}\begin{split}
&\left(\int_{\Om_+}+\int_{\Om_-}\right)(\D U \overline
V-U\overline{\D V})dxdy\\&=
-\a\int_{\Sq^1}U(0,y)\overline{\CB V(0,y)}dy.
\end{split}
\end{equation}
Indeed, substituting \eqref{1.0} into \eqref{1.gr} and
integrating by parts, we arrive at \eqref{1.gr1}.

\section{The problem of self-adjointness}\label{sad}

\vskip0.2cm In order to study  self-adjoint
realizations of $\BL_\a$ for $\a>0$, we first of all
introduce two sets,  $\CD_\a$ and $\CD_\a^\bullet$, on
which the operator is well defined. It is convenient
to do this in the terms of the expansion \eqref{1.u}.
\begin{defn}\label{dalpha}
An element $U\sim\{u_n\}\in\CM(\Om)$ lies in $\CD_\a$ if and only if
$u_n\res\R_\pm\in H^2(\R_\pm)$ for all $n\in\Z$, the conditions
\eqref{1.1x} and \eqref{1.2} are satisfied, and
\begin{equation*}
 \sum_{n\in\Z}\int_\R\bigl|-u''_n+n^2u_n\bigr|^2dx<\infty.
\end{equation*}
An element $U\in\CD_\a$ belongs to $\CD_\a^\bullet$,
if the number of non-zero terms $\{u_n\}$ in the
expansion of $U$ is finite.
\end{defn}
We denote
\[\BL_\a=-\D\res\CD_\a,\qquad
\BL^\bullet_\a=-\D\res\CD^\bullet_\a.\]

\begin{lem}\label{l1}
The operator $\BL^\bullet_\a$ is symmetric and
\[\BL_\a=(\BL^\bullet_\a)^*.\]
\end{lem}
The proof is standard and we skip it.

\begin{thm}\label{t1}
1) For $0\le\a\le 1$ the operator $\BL_\a$ is self-adjoint and,
hence, $\BL^\bullet_\a$ is essentially self-adjoint.

2) For $\a> 1$ the operator $\BL_\a$ is non-self-adjoint, and the
deficiency indices of $\BL^\bullet_\a$ are $(2,2)$.
\end{thm}

\begin{proof}
We have to check whether the equation
\begin{equation}\label{1.uu}
\BL_\a W=(\BL^\bullet_\a)^*W=\L W
\end{equation} with
$\L\neq\overline\L$ has non-zero solutions $W\in
\CD_\a$. If $W$ is such a solution, then
$W\in\CM_\L(\Om)$ and, by  \eqref{1B:Fourier}, each
component in the expansion \eqref{1.u} for $W$ can be
written as $w_n(x)=w_n e^{-|x|\sqrt{n^2-\L}}$. The
coefficients $w_n$ should satisfy
 conditions \eqref{1.2} that turn into
\begin{equation}\label{1.eq}
(n+1/2)w_{n+1}-2\a^{-1}w_n\sqrt{n^2-\L} +(n-1/2)w_{n-1}=0.
\end{equation}

The analysis of the system \eqref{1.eq} is similar to
the reasoning in  \cite{S3}, section 4, and is based
upon the classical Birkhoff -- Adams theorem, see
\cite{EL},  Theorem 8.36. The formulation of its
leading case, which we need for the study of the
operator $\BL_\a$ with $\a\neq1$, is also reproduced
in \cite{S3}. This theorem deals with one-sided
sequences ($n\in\N$ rather than $n\in\Z$ as in our
case), and we have to analyze the behaviour of $w_n$
for $n\to+\infty$ and for $n\to-\infty$ separately.

For $n\to+\infty$ we find from the theorem that
for $\a\neq1$ the equation \eqref{1.eq} has two linearly
independent solutions $\{w_n^\pm(+)\}$ such that
\begin{equation}\label{1.bap}
w_n^{\pm}(+)=(\l_+^{\pm})^n n^{-1/2}\bigl(1+O(n^{-1})\bigr),
\qquad \l_+^\pm=\a^{-1}\pm\sqrt{\a^{-2}-1}.
\end{equation}
For  $n\to-\infty$ we find in the same way that the
system has two linearly independent solutions $\{w_n^{\pm}(-)\}$
such that
\begin{equation}\label{1.bam}
w_n^{\pm}(-)=(\l_-^\pm)^n |n|^{-1/2}\bigl(1+O(|n|^{-1})\bigr),
\qquad \l_-^\pm=-\a^{-1}\pm\sqrt{\a^{-2}-1}.
\end{equation}

If $\a<1$,  we conclude from the above asymptotic
formulas that both for $n>0$ and for $n<0$ only one of
the basic solutions decays as $|n|\to\infty$. Hence,
the space of sequences $\{w_n\}$ satisfying
\eqref{1.sob} (or, equivalently, such that
$W\in\CM_\L(\Om)$) is no more than one-dimensional.
Suppose that $\{w_n\}$ is such a sequence, and apply
the following identity for solutions of recurrence
equations of the type
\[Q_{n+1}C_{n+1}+P_nC_n+Q_nC_{n-1}=0,\qquad n\in\Z,\]
with $Q_n$ real,
\begin{equation}\label{1.id}
\sum_{n=-N}^N|C_n|^2\im P_n=-Q_{N+1}\im(C_{N+1}\overline{C_N})
-Q_{-N}\im(C_{-N-1}\overline{C_{-N}}).
\end{equation}
The proof is straightforward and we skip it; cf. (4.23) and (4.24)
in \cite{S3}.

Applying \eqref{1.id} to the equation \eqref{1.eq}, we obtain
\[2\a^{-1}\sum_{n=-N}^N|w_n|^2\im\sqrt{n^2-\L}=(N+1/2)
\im(w_{N+1}\overline{w_N}+w_{-N-1}\overline{w_{-N}}).\] By
\eqref{1.bap}, \eqref{1.bam} the right-hand side vanishes as
$N\to\infty$. Since for non-real $\L$ the sign of
$\im\sqrt{n^2-\L}$ is negative if $\im\L>0$ and positive if
$\im\L<0$, we conclude that $w_n=0$ for all $n\in\Z$. It follows
that for $\a<1$ the operator $\BL_\a$ is self-adjoint.

If $\a>1$, then $|\l_\pm|=1$ and by \eqref{1.bap}, \eqref{1.bam}
any solution $\{w_n\}$ satisfies \eqref{1.sob}. This shows that for $\a>1$
the operator
$\BL_\a$ is non-self-adjoint and the deficiency indices of
$\BL^\bullet_\a$ are $(2,2)$.

\vskip0.2cm Now, let $\a=1$. Then the case ($c_1$) of the Birkhoff
-- Adams theorem applies, and the equation \eqref{1.eq} has two
linearly independent solutions of the form
\[ w_n^\pm\sim n^{\pm\sqrt{-\L}},\qquad n\to \infty\]
and similarly for $n\to -\infty$. For any non-real $\L$ only one
of such solutions may satisfy \eqref{1.sob}. Using
again the identity \eqref{1.id}, we conclude that the equation
\eqref{1.eq} has no non-zero solutions satisfying \eqref{1.sob}.
Hence, the operator $\BL_1$ is self-adjoint.
\end{proof}

\vskip0.2cm

\section{Using quadratic forms. Spectrum for  $\a\le1$}\label{smallalpha}

For small $\a$ the simplest way to study the spectrum of the
operators $\BL_\a$ is to use quadratic forms. Our argument here
follows the same line as in \cite{S3}. However, again, as in
section \ref{red}, we have to take into account that the sequence
$\{u_n\}$ is two-sided.

Integrating by parts in the expression for $(\BL_\a
U,U)$ over the semi-cylinders $\Om_\pm$, we find for
$U\in\CD^\bullet_\a$:
\begin{gather*}
(\BL_\a U,U)=\int_\Om|\nabla U|^2dxdy
 \\-\int_{\Sq^1} U'_x(0-,y)\overline{U(0,y)}dy+
\int_{\Sq^1}
U'_x(0+,y)\overline{U(0,y)}dy.\end{gather*}

Taking into account the condition \eqref{1.0}, we
obtain
\begin{gather*}
(\BL_\a U,U)-\int_\Om|\nabla U|^2dxdy
\\=i\a\int_{\Sq^1} \bigl(U'_y(0,y)\cos y+(U(0,y)\cos
y)'_y\bigr)
\overline{U(0,y)}dy\\
=i\a\int_{\Sq^1}
\bigl(U'_y(0,y)\overline{U(0,y)}-U(0,y)\overline{U'_y(0,y)}
\bigr)\cos ydy\\
=-2\a\int_{\Sq^1}\im\bigl(U'_y(0,y)\overline{U(0,y)}\bigr)\cos
ydy.\end{gather*}
 In the representation \eqref{1.u}
this turns into
\begin{equation}\label{2.1}
\bl_\a[U]:=(\BL_\a U,U)=\bl_0[U]-\a \bb[U]
\end{equation}
where
\begin{gather}
\bl_0[U]=\sum_{n\in\Z}\int_\R\bigl(|u'_n|^2+n^2|u_n|^2\bigr)dx,\label{2.2}\\
\bb[U]=\sum_{n\in\Z}
(2n-1)\re\bigl(u_n(0)\overline{u_{n-1}(0)}\bigr).\label{2.3}
\end{gather}
\vskip0.2cm Completing the set $\CD_0$ in the metric
$\bl_0[U]+\|U\|^2_{L^2(\Om)}$, we obtain a set which
we denote by $\gd$. On $\gd$ the quadratic form
$\bl_0$ is well defined and closed, and the associated
self-adjoint operator in $L^2(\Om)$ is $\BL_0$. Along
with $\gd$, we need its subspace of co-dimension one,
\begin{equation*}
  \gd'=\bigl\{U\sim\{u_n\}\in\gd:u_0(0)=0\bigr\}.
\end{equation*}

\begin{lem}\label{t2}
For any $U\in\gd'$ the following inequality is satisfied:
\begin{equation}\label{2.3b}
|\bb[U]|\le \bl_0[U]-\|u_0\|^2_{L^2(\R)}.
\end{equation}
\end{lem}
\begin{proof}
Denote by $\gd^+$ (by $\gd^-$) the subspace in $\gd$, formed by
the elements $U\sim\{u_n\}$ whose all components with $n\le0$
(with $n\ge0$) are zeroes. For $U\in\gd^\pm$  we have
$\bb[U]=\pm\bb^\pm[U]$ where
\begin{equation}\label{2.pm}
\bb^\pm[U]=\sum_{n>1}(2n-1) \re\left(u_{\pm
n}(0)\overline{u_{\pm(n-1)}(0)}\right).
\end{equation}
The estimates for $\bb^+[U]$ and for $\bb^-[U]$ are identical and
we carry them out for the `plus' sign. We derive from \eqref{2.pm}
that
\[
|\bb^+[U]|\le\sum_{n>1}
(n-1/2)\bigl(|u_n(0)|^2+|u_{n-1}(0)|^2\bigr) \le
\sum_{n\ge1}2n|u_n(0)|^2.\] Now to the $n$-th term in the last sum
we apply the elementary inequality
\[ 2\g|f(0)|^2\le\int_\R\bigl(|f'|^2+\g^2|f|^2\bigr)dx,
\qquad \forall f\in H^1(\R),\ \g>0,\]
with $\g=n$. We obtain
\begin{equation*}
|\bb^+[U]|\le \sum_{n\ge1}\int_\R
\bigl(|u'_n|^2+n^2|u_n|^2\bigr)dx.
\end{equation*}
Together with the similar inequality for $\bb^-[U]$, this yields
\eqref{2.3b}.
\end{proof}
It is not difficult to show that the factor $1$ in
front of $\bl_0[U]$ on the right-hand side of
\eqref{2.3b} cannot be improved.

\vskip0.2cm

With Lemma \ref{t2} at our disposal, it is easy to characterize
the spectral properties of the operator $\BL_\a$ for $\a<1$.

\begin{thm}\label{t3} Let $0<\a\le1$. Then

\noindent {1)}\hskip3.5cm
$\s_{ess}(\BL_\a)=\s(\BL_0)=[0,\infty)$.

\noindent {2)} The negative spectrum of $\BL_\a$ consists of
exactly one non-degenerate eigenvalue.

If $\a<1$, then also

\noindent{3)} \hskip1cm
$\s_{a.c.}(\BL_\a)=\s_{a.c.}(\BL_0)=[0,\infty),\qquad
\gm_{a.c.}(\BL_\a)=\gm_{a.c.}(\BL_0)$

\noindent(cf. \eqref{1.ac}).
\end{thm}
The proof of the statements {1)} and {3)} basically repeats the
argument in \cite{S3}, section 9, and we skip it. To justify the
statement {2)}, we first of all note that by Lemma \ref{t2}, for
$\a<1$ the quadratic form $\bl_\a$, restricted to the domain
$\gd'$, is positive definite and closed. Since $\dim\gd/\gd'=1$,
the quadratic form $\bl_\a$, considered on the whole of $\gd$, is
bounded from below and also closed. The corresponding self-adjoint
operator is $\BL_\a$. For $\a=1$, the quadratic form $\bl_\a$ is
only closable on $\gd$, and the operator $\BL_1$ corresponds to
the closure of $\bl_1$.

This reasoning shows that for $0<\a\le 1$ the number of negative
eigenvalues of $\BL_\a$ is no more than one. In order to show that
it is exactly one, it is enough to find an element $U\in\gd$ which
is such that $\bl_\a[U]<0$. To this end, we take $U\sim\{u_n\}$
with only two non-zero components $u_0,u_1$, then the desired
inequality is
\begin{equation*}
    \int_\R(|u'_0|^2+|u'_1|^2+|u_1|^2)dx<\a\re(u_1(0)\overline{u_0(0)}).
\end{equation*}
It is satisfied, for instance, if we take $u_1(x)=e^{-|x|}$ and
$u_0(x)=\vare^{-1/2}e^{-\vare|x|}$, with $\vare$ sufficiently
small.

\begin{rem}\label{2.unbb} For $\a>1$ the quadratic form $\ba_\a$ is
unbounded from below. We have to show that for any
$\a>1$ and any $M>0$ there exists an element
$U\in\gd$, such that
\begin{equation}\label{2.unbb1}
\ba_\a[U]+M\|U\|^2_{L^2(\Om)}<0.
\end{equation}
Choose a number $N\in\N$ and take $U\sim\{u_n\}$,
where $u_N(x)=e^{-|x|\sqrt{N^2+M}}$,
$u_{N-1}(x)=e^{-|x|\sqrt{(N-1)^2+M}}$ and all the
other components $u_n$ in the expansion \eqref{1.u}
are zeroes. Then
\[\begin{split}
&\int_\R\left(|u'_N|^2+(N^2+M)|u_N|^2\right)dx=2\sqrt{N^2+M},\\
&\int_\R\left(|u'_{N-1}|^2+((N-1)^2+M)|u_N|^2\right)dx=2\sqrt{(N-1)^2+M},
\end{split}\]
and
\[
\ba_\a[U]+M\|U\|^2_{L^2(\Om)}=2(\sqrt{N^2+M}+\sqrt{(N-1)^2+M})-\a(2N-1).\]
It is clear, that for any $\a>1$ the last expression
is negative, provided that $N$ is taken large enough,
and we are done.\end{rem}
\section{The case $\a>1$. Singular solutions}\label{SingSol}

In order to reach a better understanding  of
self-adjoint realizations of the operator $\BL_\a$ for
$\a>1$, we describe here the behaviour of the singular
solutions found in section \ref{sad}.

For $\a>1$ the asymptotic expressions for $w_n^\pm(\pm)$ as in
\eqref{1.bap} and \eqref{1.bam} can be re-written in a simplified
form. Indeed, set
\begin{equation*}
y(\a)=\arccos\a^{-1},
\end{equation*} then
\[\l^+_+=-\l^-_-=e^{iy(\a)},\qquad \l^+_-=-\l^-_+=e^{-iy(\a)}.\]
Therefore,
\begin{equation}\label{aux.yy}\begin{split}
&w_n^\pm(+)=e^{\pm iny(\a)}n^{-\frac12}(1+O(n^{-1})),\ n\to+\infty;\\
&w_n^\pm(-)=(-1)^ne^{\mp iny(\a)}|n|^{-\frac12}(1+O(|n|^{-1})),\
n\to-\infty.
\end{split}\end{equation}

By  \eqref{1.u}, \eqref{1B:Fourier}, and
\eqref{aux.yy}, each $L^2$-solution of the equation
\eqref{1.uu} can be represented as
\[W(x,y)=W(x,y;+)+K_0e^{-|x|\sqrt{-\L}}+W(x,y;-),\]
where $K_0$ is a constant, $W(x,y;+)$ is a certain linear
combination of the functions
\begin{equation}\label{G1.2a}\begin{split}
W^\pm(x,y;+)=&\sum\limits_{n>0}
w_n^\pm(+)n^{-\frac12}e^{iny}e^{-|x|\sqrt{n^2-\L}}\\ =
&\sum\limits_{n>0}e^{in(y\pm
y(\a))}n^{-\frac12}e^{-|x|\sqrt{n^2-\L}}
\bigl(1+O(n^{-1})\bigr),
\end{split}\end{equation}
and $W(x,y;-)$ is a linear combination of the functions
\begin{equation}\label{G1.2aa}\begin{split}
W^\pm(x,y;-)=&\sum\limits_{n<0}
w_n^\pm(-)|n|^{-\frac12}e^{iny}e^{-|x|\sqrt{n^2-\L}}\\ =
&\sum\limits_{n<0}e^{in(y\pm
y(\a)-\pi)}|n|^{-\frac12}e^{-|x|\sqrt{n^2-\L}}
\bigl(1+O(|n|^{-1})\bigr).
\end{split}\end{equation}

Note that $\sqrt{n^2-\L}=|n|+O(n^{-1})$, and hence
\[e^{-|x|\sqrt{n^2-\L}}=e^{-|x||n|}(1+|x|O(|n|^{-1})).\]
Denote by $V^\pm(x,y;\pm)$ the functions obtained by
replacing the factors $e^{-|x|\sqrt{n^2-\L}}$ by
$e^{-|x||n|}$ in each term of the sums in \eqref{G1.2a}
and \eqref{G1.2aa} and dropping the terms
$O(|n|^{-1})$. The error is a bounded function rapidly
decaying as $|x|\to\infty$. We have
 \[V^\pm(x,y;+)=\sum_{n>0}n^{-\frac12}e^{-n\left(|x|-i(y\pm y(\a))\right)}.\]

The behaviour of such sums as $|x|-i(y\pm y(\a))\to 0$
is well known. Say, it can be easily derived from the
equations (13.11) in Chapter II of the book \cite{Z}.
Denote
\[{z_\pm}^+=|x|-i(y\pm y(\a)),\qquad {z_\pm}^-=|x|-i(y\pm
y(\a)-\pi),\] then
\[  V^\pm(x,y;+)=C\left({z_\pm}^+\right)^{-\frac12}+O(1),\qquad {z_\pm}^+\to 0,\]
with an appropriate choice of the branch of the square
root,  and some constant $C$. In the same way,
\[  V^\pm(x,y;-)=C\left({z_\pm}^-\right)^{-\frac12}+O(1),\qquad {z_\pm}^-\to 0.\]

The reasoning above gives the following description of singular
solutions $W(x,y)$. These solutions depend also on the choice of $\L$,
but the leading terms of their singularities do not. For this reason
we do not reflect dependence on $\L$ in our notations.

\begin{prop}\label{SingSolu} The singular solutions $W(x,y)$ of the
equation \eqref{1.uu} have
 singularities at the points $(0,y_j)\in\Om$, were $y_j,\ j=1,2,3,4,$
are the points $\pm y(\a)$ and $\pm
y(\a)+\pi(\!\!\!\mod 2\pi)$.
  The singularity at each point is of the form
\begin{equation*}
W(x,y) \sim   C_j(|x|-iy_j)^{-\frac12}+O(1).
\end{equation*}
\end{prop}

In order to explain the role of these four singular points, let us
check the Shapiro -- Lopatinsky criterion for the ellipticity of
the boundary-value problem $-\D U=F$ under the conditions
\eqref{1.-1} and \eqref{1.0}. In our case this criterion
determines the point $y\in\Sq^1$ as regular if and only if the
problem
\begin{equation*}\label{G7.Lopat}
  -\phi''(t)+\phi(t)=0,\; t\ne0,\qquad \phi'(0+)-\phi'(0-)=\pm 2\a
  \cos y\phi(0)
\end{equation*}
has only trivial bounded continuous solutions on the
line $t\in(-\infty,\infty)$. This requirement is
violated exactly at the points $y=y_j,\ j=1,2,3,4,$
where $\a|\cos y|=1$, the solution being
$\phi(t)=e^{-|t|}$. On the other hand, for
$\a\in[0,1)$ the Shapiro -- Lopatinsky condition is
satisfied at all transition points. Therefore every
weak solution of the equation $-\Delta U=\L U$
satisfying \eqref{1.0} belongs to $H^2$ in both
half-cylinders $\Omega_\pm$, so it is non-singular,
which explains the self-adjointness.

\section{The case $\a>1$. Spectral properties}\label{BigAlfaspec}

For $\a>1$, the main technical difficulty stems from
the fact that  Definition~\ref{dalpha} does not
describe the class $\CD_\a$ in the terms of standard
function spaces on $\Om$. For this reason, our
argument here is rather lengthy.

Let us fix some self-adjoint extension $\hat\BL_\a$ of
the operator $\BL_\a^\bullet$. The spectral properties
discussed in this section do not depend on the choice
of the extension.

We start by establishing a formula for the difference of
resolvents of the operators $\hat\BL_\a$ and $\BL_0$. The method
for finding this kind of expressions is widely used and was
proposed by  Birman in \cite{Bir}. Let first $\L$ be a
non-real number. It belongs to the resolvent sets of both
operators $\hat\BL_\a$ and $\BL_0$, and we denote by $\hat\BR_\a$,
$\BR_0$ the corresponding resolvents.

Take some $F,G\in L^2(\Omega)$, and consider the sesqui-linear
form
\begin{equation}\label{G7.form1}
\br[F,G]=((\hat\BR_\a-\BR_0)F,G)=(\hat\BR_\a F,G)-(F,\BR_0^*G).
\end{equation}
Denote
 \[\hat\BR_\a F=U,\qquad \BR_0^*G=V,\]
then $U\in\CD_\a$ and $V\in H^2(\Omega)$. Thus the quadratic form
\eqref{G7.form1} can be re-written as
\begin{equation*}
(U, (\BL_0-\overline{\L}) V)-((\hat\BL_\a-\L) U, V)=
\left(\int_{\Om_+}+\int_{\Om_-}\right) (\D
U\overline{V}-U\overline{\D V})dxdy.
\end{equation*}

Applying \eqref{1.gr1}, we arrive at
\begin{equation*}
\br{[F,G]}=\a\int_{\Sq^1} U(0,y)  \overline{\CB V(0,y)}dy,
\end{equation*}
where $\CB$ is the operator \eqref{1.b}.
Hence, the latter equality gives the representation of the operator
$\hat\BR_\a-\BR_0$ as

\begin{equation}\label{G7.operator4}
\hat\BR_\a-\BR_0= 2\a S^*T , \qquad T=\G\hat\BR_\a,\ S=\CB\G\BR_0^*,
\end{equation}
where $\G$ stands for the operator of restriction of  functions
on $\Om$ to the circle $x=0$. The operator $T$ is bounded from
$L^2(\Omega)$ to $H^{-\frac12}(\Sq^1)$, and $S$ is bounded from
$L^2(\Omega)$ to $H^{\frac12}(\Sq^1)$, so that $S^*$ is bounded
from $H^{-\frac12}(\Sq^1)$ to $L^2(\Omega)$.

Our next step is to derive a pseudo-differential
equation for the distribution $w=\G W$, where
\begin{equation}\label{5.xx}
W=U-V_1:=\hat{\BR}_\a F-\BR_0 F,\qquad F\in L^2(\Om).
\end{equation}
Evidently, $W\in\CM_\L(\Om)$ and thus $w\in
H^{-\frac12}(\Sq^1)$. Below we denote by $\CA$ the
operator
 $-\frac{d^2}{dy^2}$ in $L^2(\Sq^1)$, extended to
 distributions on $\Sq^1$. It follows from the representation
\eqref{1B:Fourier} that
\[W'_x(0+,y)-W'_x(0-,y)=-2\sum_nw_n\sqrt{n^2-\L}e^{iny}=
-2(\CA-\L)^{\frac12}w(y).\] Now, taking into account
the transmission conditions for $U$ and for $V_1$, we
find that
\[W'_x(0+,y)-W'_x(0-,y)-\a \CB W(0,y)=\a \CB V_1(0,y),\]
or
\begin{equation}\label{5.eq}
\left(2(\CA-\L)^{\frac12}+\a \CB\right)w=-\a \CB\G
\BR_0F\in H^{\frac12}(\Sq^1).
\end{equation}

\vskip0.2cm

The operator $\hat\BR_\a-\BR_0$ is, of
course, bounded.
 We are going to show that, actually, it is compact.
 The proof is based  upon the fact that the operator $T$
 in \eqref{G7.operator4} acts from $L^2(\Omega)$ not only into
$H^{-\frac12}(\Sq^1)$ but into a smaller
 space, $H^{-\e}(\Sq^1)$, for any $\e>0$.
 To show  this, we need an  {\it a priori} estimate for the equation \eqref{5.eq}.
This equation is elliptic for $\a<1$, but for $\a>1$, which is the case we
are dealing with, it is degenerate, so some more effort is needed.

 \begin{lem}\label{G7.lemma1}
For any $\e>0$  there exist constants $C,C'$  such that
  for any  $w\in H^{-\frac12}(\Sq^1)$
 \begin{equation}\label{G7.apriori1}
   \|w\|_{H^{-\e}(\Sq^1)}\le
   C\|2(\CA-\L)^{\frac12}w+\a \CB w\|_{H^{\frac12}(\Sq^1)}+
C'\|w\|_{H^{-\frac12}(\Sq^1)},
\end{equation}
provided that the first term on the right-hand side of
\eqref{G7.apriori1} is finite.
\end{lem}

\begin{proof}
Denote by $P_\pm$ the Riesz projections,
 \[   P_+f=\pi^{-1}\sum_{k\ge 0} (f,e^{iky})e^{iky},\qquad
    P_-f=\pi^{-1}\sum_{k< 0} (f,e^{iky})e^{iky}.\]
Here the sums are understood in the sense of distributions; in
particular, if $f\in H^{s}(\Sq^1), \; s\in \R$, both series
converge in $H^{s}(\Sq^1)$.

The operators $P_\pm$ differ by smoothing operators from
pseudodifferential operators on the circle with symbols
\[ p_{+}(y,\y)=\begin{cases} 1&{\text{if}}\ \y>0,\\ 0&{\text{if}}\
\y<0;\end{cases}\qquad p_-(y,\y)=1-p_{+}(y,\y),\]
 see the
discussion in \cite{Agr} about the Fourier series representation
of pseudodifferential operators on the circle.

 For $w\in H^s(\Sq^1)$ we  denote by $w_\pm$
the distributions $w_\pm=P_\pm w$.  The operator
$(\CA-\L)^{\frac12}$ is, up to a smoothing term,  the
pseudodifferential operator with symbol
$(\y^2-\L)^{\frac12}=|\y|+O(|\y|^{-1})$. As it follows
from the composition formulas for pseudodifferential
operators in dimension one, the operators in
\eqref{G7.apriori1} commute or almost commute with
$P_\pm$:
\[(\CA-\L)^{\frac12}P_\pm=P_\pm(\CA-\L)^{\frac12},\qquad \CB P_\pm=P_\pm \CB
+K,\] with $K$ being a smoothing operator.  Thus, up to en error
being an operator of order $-1$, the operator $(\CA-\L)^{\frac12}$
acts on the components $w_\pm$ as the differentiation, with proper
coefficients:
$$\|(\CA-\L)^{\frac12}w_\pm \mp iw_\pm'\|_{H^{\frac12}(\Sq^1)}\le
C\|w_\pm\|_{H^{-\frac12}(\Sq^1)}.$$
Therefore, \eqref{G7.apriori1} will
follow as soon as we prove that
\begin{equation}\label{G7.separated}
    \|w_\pm\|_{H^{-\e}(\Sq^1)}\le
     C\|\pm w_\pm'+\frac{\a}{2i}
     \CB w_\pm\|_{H^{\frac12}(\Sq^1)}+C'\|w_\pm\|_{H^{-\frac12}(\Sq^1)}.
\end{equation}

The estimate \eqref{G7.separated}, even  with $-\e$
replaced by $3/2$ on the left-hand side, would follow
automatically, if the operators $\pm 2i
\partial_y +\a \CB$ were elliptic for both signs $\pm$.
This is the case for $|\a|<1$. However for $|\a|\ge1$
these operators have points of degeneracy of
ellipticity, i.e. the points where the principal
symbols $(\pm 1+ \a \cos y) \eta$ vanish. Note that
these are exactly the points where the singularities
of the singular solutions are located, see section
\ref{SingSol}. For such degenerate operators
considering  the principal symbol is not sufficient
for getting {\it a priori} estimates, so the influence
of lower order terms in $\CB$ must be taken into
account.

We concentrate on the case of the 'minus' sign in
\eqref{G7.separated}. Let us denote $h(y)=\a\cos y$ and set
\[u=-w_-'+\frac{\a}{2i}\CB w_-=(h(y)-1)w_-'+\frac12 h'(y)w_-.\]
 We also set $g=(h(y)-1)^\frac12 w_-$, with a properly chosen branch
of the square root. Note that $g'=(h(y)-1)^{-\frac12}u$. Our next
task is to derive an estimate of $g$ in the terms of $u$, assuming that $u \in H^{\frac12}(\Sq^1)$.
\vskip0.2cm
The latter assumption on $u$ implies that the function
$(h(y)-1)^{-\frac12}u$ belongs to the space $H^{-\d}(\Sq^1)$ for an
arbitrarily small $\d>0$, say $\d<1/2$. To justify the above
statement, we must show that
\begin{equation}\label{5.functional}
    \left|\int_{\Sq^1}(h(y)-1)^{-\frac12}u(y)\z(y)dy\right|\le C
    \|u\|_{H^{\frac12}(\Sq^1)}\|\z\|_{H^\d(\Sq^1)},\qquad \forall\z\in H^\d(\Sq^1).
\end{equation}
But this follows from the H\"older inequality, since
$|h-1|^{-\frac12}\in L^r(\Sq^1)$ for any $r<2$, and by the embedding
theorem $u\in L^q(\Sq^1)$ for any $q<\infty$ and
$\z\in L^{\frac{2}{1-2\d}}(\Sq^1)$.

It follows from \eqref{5.functional} that
$$\| g'\|_{H^{-\d}(\Sq^1)}=\|(h(y)-1)^{-\frac12}u\|_{H^{-\d}(\Sq^1)}\le
 C\|u\|_{H^{\frac12}(\Sq^1)}.$$
 Therefore, the function $g=(h(y)-1)^\frac12 w_-$ lies in $H^{1-\d}(\Sq^1)$ and
 satisfies the estimate
 $$\|g\|_{H^{1-\d}(\Sq^1)}\le C \|u\|_{H^{\frac12}(\Sq^1)}+
C'\|g\|_{H^{-N}(\Sq^1)},$$
with $N$ being arbitrarily large.

By the definition of $g$, we have
$w_-=(h(y)-1)^{-\frac12}g$. An estimate, similar to
\eqref{5.functional} (even a simpler one, since $g\in
L^\infty(\Sq^1)$), shows that $w_-$  belongs to
$H^{-\e}(\Sq^1)$, with the required estimate.
\end{proof}

 The estimate, just proved, enables us to establish
 the
 compactness of the difference of resolvents
 $\hat\BR_\a-\BR_0$ and of several related operators and
 prove spectral estimates.

  \begin{prop}\label{G7.Prop} The operator
  $\hat\BR_\a-\BR_0$ is compact, moreover for its singular
  numbers $s_n(\hat{\BR}_\a-\BR_0)$ the estimate
  \begin{equation}\label{G7.Prop.Differ}
s_n(\hat\BR_\a-\BR_0)=O(n^{-\frac12+\e})
\end{equation}
holds for any $\e>0$. Further on,
\begin{equation}\label{G7.Prop.DifferR}
s_n((\hat\BR_\a-\BR_0)\BR_0)=O(n^{-\frac52+\e}),\qquad
s_n(\BR_0(\hat\BR_\a-\BR_0))=O(n^{-\frac52+\e}).
\end{equation}
\end{prop}

\begin{proof}
It follows from   the factorization
\eqref{G7.operator4} that
\[s_n(\hat\BR_\a-\BR_0)\le Cs_n(T),\]
where  we have to consider the operator $T$ as
acting from $L^2(\Omega)$ to $H^{-\frac12}(\Sq^1)$. For $F\in
L^2(\Omega)$ we define the function $W$ as in \eqref{5.xx}
and take $w=W(0,\cdot)$, then
 $$TF=\G\hat{\BR}_\a F=w+\G\BR_0F.$$
The operator $\G\BR_0$ acts from $L^2(\Omega)$ to $H^{\frac32}(\Sq^1)$,
and the distribution $w$ satisfies the equation \eqref{5.eq}, whose right-hand
side belongs to $H^{\frac12}(\Sq^1)$. Lemma \ref{G7.lemma1} applies
and gives $w\in H^{-\e}(\Sq^1)$. It follows that the operator $T$ is bounded
as acting from $L^2(\Omega)$ to $H^{-\e}(\Sq^1)$, and therefore, the
singular numbers of the operator $T:L^2(\Omega)\to H^{-\frac12}(\Sq^1)$
are controlled by those of the embedding $H^{-\e}(\Sq^1)\to H^{-\frac12}(\Sq^1)$.
The latter are of the order $O(n^{-\frac12+\e})$, whence the
required estimate \eqref{G7.Prop.Differ}.

\vskip0.2cm

 Further on, we factorize the operator $\BR_0(\hat\BR_\a-\BR_0)$ as
\[\BR_0(\hat\BR_\a-\BR_0)=2\a\BR_0 S^*T.\]
Since we already know the singular numbers estimate
for the operator $T:L^2(\Om)\to H^{-\frac12}(\Sq^1)$,
it is sufficient for us to consider the operator
$\BR_0 S^*$ as acting between the spaces
$H^{-\frac12}(\Sq^1)$ and $L^2(\Om)$. It is more
convenient to deal with the adjoint operator
\begin{equation}\nonumber
S\BR_0^*=\CB\G{(\BR_0^*)}^2: L^2(\Om)\to H^{\frac12}(\Sq^1).
\end{equation}
This operator is bounded as acting from $L^2(\Om)$ to
$H^{\frac52}(\Sq^1)$. Hence, the singular numbers of
the same operator but considered as acting between the
spaces $L^2(\Om)$ and $H^{\frac12}(\Sq^1)$ are
controlled by those of the embedding operator
$H^{\frac52}(\Sq^1)\to H^{\frac12}(\Sq^1)$. The latter
are of the order $O(n^{-2})$. This, together with the
estimate for $T$, proves the second estimate in
\eqref{G7.Prop.DifferR}. The first estimate in
\eqref{G7.Prop.DifferR} follows from the second one by
passing to adjoint operators.
\end{proof}

Now we arrive at our main result on the spectrum of the operator
$\hat{\BL}_\a$, $\a>1$.
\begin{thm}\label{G7:theorem} For $\a>1$ the spectrum
of the operator $\hat{\BL}_\a$ consists of the
essential spectrum filling the semi-axis $\l\ge0$ and
the eigenvalues below the point $0$. The set of
eigenvalues below the essential spectrum is unbounded
from below, may have only $0$ and $-\infty$  as limit
points, and for the counting function $n(t)=\#\{\l\in
\s_{\mathrm{disc}}(\hat{\BL}_\a),\l\in(-t,-t_0)\}$,
with any fixed $t_0>0$, the estimate holds
\begin{equation}\label{G7:counting function}
n(t)=O(t^{2+\e_1}), \mathrm{\;for\; any\;} \e_1>0.
\end{equation}
The absolutely continuous spectrum of $\hat{\BL}_\a$
fills the half-line $\l\ge0$ and its multiplicity
function coincides with that of $\BL_0$.
\end{thm}

\begin{rem} The estimate \eqref{G7:counting function} is rather rough.
The authors believe that a more detailed analysis,
based upon a further study of the degenerate equation
\eqref{5.eq}, would show that the counting function
has the asymptotics $n(t)\sim C t^{\frac12}$ as $t\to
\infty$. Moreover, we think that the negative
eigenvalues do not have 0 as their limit point.
\end{rem}

\begin{proof} First, we note that
due to Weyl theorem, the essential spectrum of the operators
$\hat\BR_\a$ and
 $\BR_0$ is the same, therefore
  the essential spectrum of $\hat\BL_\a$ coincides with that of  $\BL_0$,
so it is
 the half-line $[0,\infty)$. Thus, the
  spectrum of $\hat\BL_\a$ below $0$ may only
consist of eigenvalues with possible accumulation points only at
$0$ and $-\infty$. The latter point must be
 an accumulation point for eigenvalues since the operator $\hat\BL_\a$ is not
 semi-bounded from below, see Remark~\ref{2.unbb}.
 The discreteness of the negative spectrum implies that there are real
 regular points of the operator $\hat{\BL}_\a$, these are all points below $0$,
  which are not eigenvalues. We fix such  regular
  $\L<0$ and consider the resolvents $\hat\BR_\a,\BR_0$ at
  this point.
  Then the above
construction of the
 operator $\hat\BR_\a-\BR_0$ and the estimate
 for its eigenvalues can be repeated, this time for the chosen real $\L$.
The spectrum of $\BR_0$ coincides with the interval
$[0,-\L^{-1}]$, and
\begin{equation*}
\hat\BR_\a=\BR_0+(\hat\BR_\a-\BR_0).
\end{equation*}
The operator $\BR_0$ is non-negative, therefore, for any $\mu<0$
the number of eigenvalues of $\hat\BR_\a$ (counting
multiplicities)
 in $(-\infty,\mu)$ is not greater than
 the number of eigenvalues of $\hat\BR_\a-\BR_0$ in the same interval. The
  latter quantity is estimated by means
 of the eigenvalue bound \eqref{G7.Prop.Differ}, which
 under an appropriate choice of $\e=\e(\e_1)$
 leads to \eqref{G7:counting function}, with $t_0=-\L$ and $t=-(\mu^{-1}+\L)$.

 In order to justify the statement on the absolute
 continuous spectrum, let us consider the difference
 $\hat\BR_\a^3-\BR_0^3.$ We have
\begin{gather*}
\hat\BR_\a^3-\BR_0^3= (\hat\BR_\a-\BR_0)^3+\hat\BR_\a\BR_0(\hat\BR_\a-\BR_0)\\
+\hat\BR_\a(\hat\BR_\a-\BR_0)\BR_0+\BR_0(\hat\BR_\a-\BR_0)^2+\BR_0^2
(\hat\BR_\a-\BR_0),
\end{gather*}
and due to the estimates   \eqref{G7.Prop.Differ} and
\eqref{G7.Prop.DifferR} each term is trace class. By
Kato's theorem, the absolute continuous parts of
operators $\hat{\BL}_\a$ and ${\BL}_0$ are unitary
equivalent.
\end{proof}

\section{An alternative model}\label{alt}
Here we briefly describe an alternative model, where a
slight change in the setting leads to some major
changes in the spectral behaviour. The family $\BM_\a$
of differential operators acts on the strip
$\Om'=\R\times(0,\pi)$ and is generated by the
Laplacian $-\D U=-U''_{x^2}-U''_{y^2}$, the Dirichlet
condition $U(x,0)=U(x,\pi)=0$, and two additional
conditions at $x=0$:
\begin{equation*}
\begin{split}
 U(0+,y)&=U(0-,y)\ (=U(0,y)),\\
U'_x(0+,y)-U'_x(0-,y)&=-i\a\left(U'_y(0,y)\sin y +(U(0,y)\sin
y)'_y\right),\end{split}
\end{equation*}
cf. \eqref{1.-1}, \eqref{1.0}.
 \vskip0.2cm

The Fourier expansion for this case has the form
\[U=\sum_{n=1}^\infty u_n(x)\varf_n (y),\qquad \varf_n(y)=\sqrt{\frac{2}{\pi}}\sin ny\]
(in short, $U\sim\{u_n\}$). The equation and the boundary and
transmission conditions reduce to an infinite system of ordinary
differential operators on the real axis, coupled by the
conditions at $x=0$:
\begin{equation}\nonumber
-\D U\sim\{-u_n''+n^2u_n\},\qquad n\in\N;
\end{equation}
each $u_n$ is continuous at $x=0$;
\begin{equation}\nonumber
u_n'(0+)-u_n'(0-)=i\a\bigl((n+1/2)u_{n+1}(0)-(n-1/2)u_{n-1}(0)\bigr),
\end{equation} with $u_0$ taken to be
identically zero.

For $\a=0$ the system decouples, and we get an analogue of
\eqref{1.l0}, but this time with summation over $n\in\N$. From
here we derive that the spectrum $\s(\BM_0)$ is absolutely
continuous, fills the half-line $[1,\infty)$, and its multiplicity
function is given by
\begin{equation}\label{1A.ac}
    \gm_{a.c.}(\l;\BM_0)=2[\l], \qquad\forall \l\ge 1.
\end{equation}

It is these two differences with $\BL_0$, the sequence of $u_n$
being one-sided and the spectrum of the unperturbed problem
starting at $1$ rather than  at $0$, that lead to the changes in
the spectral properties of the perturbed operator.

\vskip0.2cm The  study of the self-adjointness of
$\BM_\a$ for $\a>0$ follows the same line as for the
operators $\BL_\a$ in section \ref{sad}. It turns out
that the operator $\BM_\a$, considered on the natural
domain (cf. Definition \ref{dalpha}), is self-adjoint
for $\a\le 1$. If $\a>1$, the operator has a
one-parameter family $\hat\BM_\a$ of self-adjoint
realizations. The singular solutions, which define
these realizations by v.Neumann's scheme, have two
singular points $(0,y^{\pm})$, with singularities of
the order $C(|x|+i(y-y^{\pm}))^{-\frac12}$. The points
$y^\pm$ are the solutions of the equation $\a\sin
y=1$, these are exactly the points where the Shapiro
-- Lopatinsky condition is violated.

Similarly to the cylinder case, the spectral analysis of the
operator $\BM_\a$ for $0<\a< 1$ is based upon considering the
quadratic forms. The quadratic form for $\BM_\a$ is
\begin{equation}\nonumber
\bm_\a[U]:=(\BM_\a U,U)=\bm_0[U]-\a \bb[U]
\end{equation}
where
\begin{gather*}
\bm_0[U]=\sum_{n\in\N}\int_\R\bigl(|u'_n|^2+n^2|u_n|^2\bigr)dx,\\
\bb[U]=\sum_{n\ge
2}(2n-1)\im\bigl(u_n(0)\overline{u_{n-1}(0)}\bigr),
\end{gather*}
cf. \eqref{2.1}, \eqref{2.2} and \eqref{2.3}.
 The quadratic form $\bm_0$ is positive definite and
closed on its natural domain which we again denote by $\gd$. The
associated self-adjoint operator in $L^2(\Om')$ is $\BM_0$. The
inequality
\begin{equation}\label{2A.3a}
|\bb[U]|\le \bm_0[U], \qquad U\in \gd,
\end{equation}
is checked in the same way as \eqref{2.3b}, and this
time no second term as in \eqref{2.3b}  appears. The
constant factor $1$ in the estimate \eqref{2A.3a} is
sharp. Hence, for $\a<1$ the quadratic form $\bm_\a$
is positive definite and closed on $\gd$. The
corresponding self-adjoint operator in $L^2(\Om')$ is
$\BM_\a$. It is not difficult to show that for $\a>1$
the quadratic form $\bm_\a$ is unbounded from below.

\bigskip

We pass now to the description of the spectrum of
$\BM_\a$. It is here where the differences with
$\BL_\a$ manifest themselves, cf. \thmref{t3}.

\begin{thm}\label{t3A}
Let $0<\a<1$. Then

\noindent 1)\hskip3.5cm $\s_{ess}(\BM_\a)=\s(\BM_0)=[1,\infty)$.

\noindent 2)\hskip1.3cm
$\s_{a.c.}(\BM_\a)=\s_{a.c.}(\BM_0)=[1,\infty),\qquad
\gm_{a.c.}(\BM_\a)=\gm_{a.c.}(\BM_0)$

\noindent (cf. \eqref{1A.ac}).

\noindent 3) The spectrum of $\BM_\a$ below the threshold $\l_0=1$
is finite.
\end{thm}
We skip the proof which basically repeats the argument
in \cite{S3}, section 9. Note that one can also prove
that for the pairs $\BM_\a, \BM_0$ and $\BM_0, \BM_\a$
there exist complete isometric wave operators.

The quadratic form $\bm_1\res\gd$ is non-negative and closable, it
generates the
operator $\BM_1$. It is possible to show that its essential spectrum
is the half-line $[0,\infty)$.

\vskip0.2cm

The  analysis of the discrete spectrum of $\BM_\a$ for
$\a\in (0,1)$ is based upon a version of
Birman-Schwinger principle found in \cite{S3}. Before
giving its formulation, let us recall the following
well-known notations. Given a real number $\l$ and
self-adjoint operator $Q$, whose spectrum on
$(-\infty,\l)$ is discrete, we write $N_-(\l;Q)$ for
the number of the eigenvalues $\l_n(Q)<\l$, counted
according to their multiplicities. We also write
$N_+(\l;Q)=N_-(-\l;-Q)$.

It turns out that within an error which is no greater than $1$,
the number $N_-(1;\BM_\a)$ coincides with $N_+(\a^{-1};\BJ)$,
where $\BJ$ is a certain infinite Jacobi matrix:
\begin{equation}\label{4A.6}
0\le N_-(1;\BM_\a)-N_+(\a^{-1};\BJ)\le1.
\end{equation}
The reasoning is the same as in \cite{S3}, however the Jacobi
matrix $\BJ$ turns out to be different: it is the zero-diagonal
Jacobi matrix, with the non-diagonal entries given by
\[ 2j_{n,n-1}=2j_{n-1,n}=\frac{n-1/2}{(n^2-1)^{1/4}(n^2-2n)^{1/4}}.\]
Since $j_{n,n-1}\to 1/2$ as $n\to\infty$, the matrix $\BJ$ has the
absolutely continuous spectrum filling the segment $[-1,1]$ and
the spectrum outside this segment is discrete. Note that
$\a\in(0,1)$ is equivalent to $\a^{-1}>1$, so that both terms in
\eqref{4A.6} are finite.

In order to estimate $N_+(\mu;\BJ)$, $\mu=\a^{-1}$, we use the
asymptotics of $j_{n,n-1}$:
\begin{equation}\label{AsJac}
j_{n,n-1}\sim \frac12 +\frac12 n^{-2} +o(n^{-2}), \; n\to\infty.
\end{equation}
 Using the results of Geronimo
\cite{Ger1}, \cite{Ger2}, combined with some standard variational
tools, one can show that $N_+(\mu;\BJ)$ can be estimated from
below and from above by $|\log(\mu-1)|$, with different constants.
Thus the number of eigenvalues of $\BM_\a$ in $(0,1)$ grows
logarithmically as $\a\nearrow 1$. We believe that actually a
logarithmical asymptotics for the eigenvalues holds.

When $\a$ becomes larger than $1$, the phase transition occurs,
similar to the cylinder case. Each self-adjoint realization
$\hat{\BM}_\a$ of the operator $\BM_\a$ is unbounded from below,
with the spectrum below the point $1$ being discrete. The
absolutely continuous spectrum is still the half-line
$[1,\infty)$, with the same multiplicity function as for $\BM_0$.
All these properties are proved using the methods exposed in
section~\ref{BigAlfaspec}. Some additional technical complications
are  caused by the fact that now we should prove estimates of the
type \eqref{G7.apriori1} for the operators on an interval
$(0,\pi)$, rather than on the circle $\Sq^1$ which is a manifold
without boundary. But these complications can be overcome.

\section{Acknowledgements}
The work on the paper was started in April of 2005, when G.
Rozenblum visited the Weizmann Institute of Science. G.R.
expresses his gratitude to the Institute for its hospitality and
financial support.

The authors are also grateful to Y. Kannai for a very useful
discussion.

\bibliographystyle{amsalpha}

\end{document}